\documentclass[a4paper,12pt,reqno]{amsart}
\usepackage[a4paper,margin=1.9cm,top=2.5cm,bottom=2.5cm,centering,vcentering]{geometry}
\usepackage{graphics}
\usepackage{graphicx}
\usepackage{amssymb}
\usepackage{pstricks}\usepackage{tikz}
\usepackage[colorlinks=true]{hyperref}

\allowdisplaybreaks[1]

\newcounter{c-prop}
\newcounter{c-conn}
\newcounter{c-disc}
\newcounter{c-gr}

\newtheorem{theorem}{Theorem}

\DeclareMathOperator{\conv}{\text{conv}}
\newcommand{\e}{\tilde{e}}
\newcommand{\A}{\mathcal{A}}
\newcommand{\E}{\mathcal{E}}
\newcommand{\G}{\mathcal{G}}
\newcommand{\Gc}{\mathcal{G}_{\text{C}}}
\renewcommand{\P}{\mathcal{P}}
\newcommand{\PP}{\widetilde{\mathcal{P}}}
\newcommand{\R}{\mathbb{R}}
\newcommand{\Nsing}[1]{\#(\text{single\:edges},#1)}
\newcommand{\Ndoub}[1]{\#(\text{edge\:pairs},#1)}
\newcommand{\Nloops}[1]{\#(\text{loops},#1)}
\newcommand{\Gtree}{\G_{\text{trees}}}
\newcommand{\Gltree}{\G_{\text{looped\,trees}}}
\newcommand{\Getree}{\G_{\text{enh.\,trees}}}
\newcommand{\Gqtree}{\G_{\text{quasitrees}}}

\DeclareMathOperator{\ehr}{Ehr}
\DeclareMathOperator{\vol}{vol}

\title[Ehrhart polynomials of partial permutohedra]{Ehrhart polynomials of partial permutohedra}
\author[R.~E.~Behrend]{Roger E.~Behrend}
\address{R.~E.~Behrend, School of Mathematics, Cardiff University, Cardiff,
CF24 4AG, UK}
\email{behrendr@cardiff.ac.uk}
\begin{document}
\begin{abstract}
For positive integers $m$ and $n$, the partial permutohedron~$\mathcal{P}(m,n)$ is a certain integral polytope in $\mathbb{R}^m$,
which can be defined as the convex hull of the vectors
from $\{0,1,\ldots,n\}^m$ whose nonzero entries are distinct.
For $n=m-1$, $\mathcal{P}(m,m-1)$ is (after translation by $(1,\ldots,1)$) the polytope~$P_m$ of parking functions of length~$m$,
and for $n\ge m$, $\mathcal{P}(m,n)$ is combinatorially equivalent to an $m$-stellohedron.
The main result of this paper is an explicit expression for the Ehrhart polynomial of~$\mathcal{P}(m,n)$ for any $m$ and $n$ with $n\ge m-1$.
The result confirms the validity of a conjecture for this Ehrhart polynomial in arXiv:2207.14253~\cite{BehCasChaDiaEscHarInk22}, and the $n=m-1$ case also
answers a question of Stanley regarding the number of integer points in $P_m$.
The proof of the result involves transforming $\mathcal{P}(m,n)$ to a unimodularly equivalent polytope in $\mathbb{R}^{m+1}$,
obtaining a decomposition of this lifted version of $\mathcal{P}(m,n)$ with $n\ge m-1$ as a Minkowski sum of dilated coordinate simplices,
applying a result of Postnikov for the number of integer points in generalized permutohedra of this form,
observing that this gives an expression for the Ehrhart polynomial of $\mathcal{P}(m,n)$ with $n\ge m-1$ as an edge-weighted sum
over graphs (with loops and multiple edges permitted) on~$m$ labelled vertices in which each connected component contains at most one cycle, and
then applying standard techniques for the enumeration of such graphs.
\end{abstract}
\maketitle

\section{Introduction}\label{intro}
For positive integers~$m$ and~$n$, the partial permutohedron~$\P(m,n)$ 
can be defined as the convex hull 
of the vectors from $\{0,1,\ldots,n\}^m$ whose nonzero entries are distinct.
This class of integral polytopes in $\R^m$ was introduced and studied by Heuer and Striker~\cite{HeuStr21,HeuStr22},
and has been studied further by Behrend, Castillo, Chavez, Diaz-Lopez, Escobar, Harris and Insko~\cite{BehCasChaDiaEscHarInk22,BehCasChaDiaEscHarInk23},
Black and Sanyal~\cite{BlaSan22} (in the context of monotone path polytopes of polymatroids), and Hanada, Lentfer and Vindas-Meléndez~\cite{HanLenVin23}
(in the context of generalized parking function polytopes).

As examples of partial permutohedra, $\P(3,n)$ can be depicted as
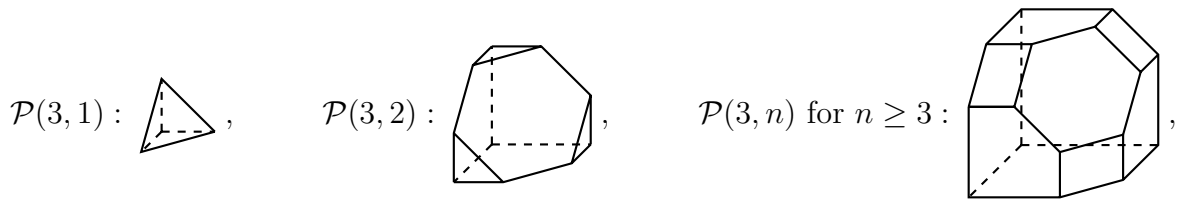
\begin{equation*}\P(3,1):\;\raisebox{-4mm}{\begin{tikzpicture}[x={(.7,0,0)},y={(0,.7,0)},z={(0,0,.7)}]\draw[thick](1,0,0)--(0,1,0)--(0,0,1)--(1,0,0);\draw[dashed,thick]	(0,0,0)--(1,0,0)(0,0,0)--(0,1,0)(0,0,0)--(0,0,1);\end{tikzpicture}}\;,\hspace{11mm}
\P(3,2):\;\raisebox{-8mm}{\begin{tikzpicture}[x={(.65,0,0)},y={(0,.65,0)},z={(0,0,.65)}]
\draw[thick] (2,0,0)--(2,1,0)--(1,2,0)--(0,2,0)(2,0,0)--(2,0,1)--(1,0,2)--(0,0,2)(0,0,2)--(0,1,2)--(0,2,1)--(0,2,0)(0,1,2)--(1,0,2)(2,0,1)--(2,1,0)(0,2,1)--(1,2,0);
\draw[dashed,thick](0,0,0)--(2,0,0)(0,0,0)--(0,2,0)(0,0,0)--(0,0,2);\end{tikzpicture}}\;,\hspace{11mm}
\P(3,n)\text{ for }n\ge3:\;\raisebox{-10mm}{\begin{tikzpicture}[x={(.6,0,0)},y={(0,.6,0)},z={(0,0,.6)}]
\draw[thick](1,2,3)--(2,1,3)--(3,1,2)--(3,2,1)--(2,3,1)--(1,3,2)--(1,2,3)(1,2,3)--(0,2,3)(2,1,3)--(2,0,3)(3,1,2)--(3,0,2)(3,2,1)--(3,2,0)(2,3,1)--(2,3,0)(1,3,2)--(0,3,2)(1,2,3)--(0,2,3)(0,2,3)--(0,0,3)--(2,0,3)--(3,0,2)--(3,0,0)--(3,2,0)--(2,3,0)--(0,3,0)--(0,3,2)--(0,2,3);
\draw[dashed,thick] (0,0,0)--(3,0,0)(0,0,0)--(0,3,0)(0,0,0)--(0,0,3);\end{tikzpicture}}\;,\end{equation*}
where (as will be discussed in~(viii) 
in Section~\ref{prop}) all cases with $n\ge3$ are combinatorially equivalent.

The main aim of this paper is to prove the following result, which was conjectured in~\cite[Conj.~6.5]{BehCasChaDiaEscHarInk22}.

\begin{theorem}\label{EhrTh}
For any positive integers~$m$ and~$n$ with $n\ge m-1$, the Ehrhart polynomial $\ehr_{\P(m,n)}(t)$ of~$\P(m,n)$ is explicitly
\begin{multline}\label{EhrEq}
\ehr_{\P(m,n)}(t)\big|_{n\ge m-1}\\
=\frac{1}{2^m}\sum_{i=0}^{\lfloor m/2\rfloor}\sum_{j=2i}^m(-1)^{i+1}\,\binom{m}{m-j,\,j-2i,\,i,\,i}\,i!\:(2j-4i-3)!!\:t^{j-i}\,(2nt+t+2)^{m-j}.\end{multline}
\end{theorem}
In~\eqref{EhrEq}, standard notation is used for multinomial coefficients and double factorials, i.e., $\binom{m}{m-j,\,j-2i,\,i,\,i}=m!/\bigl((m-j)!\,(j-2i)!\,(i!)^2\bigr)$ and $(2j-4i-3)!!=-\prod_{k=1}^{j-2i}(2k-3)$.
For general information regarding Ehrhart polynomials of integral polytopes, see for example~\cite[Ch.~3]{BecRob15} or~\cite[Sec.~4.6.2]{Sta12}.

\smallskip
The proof of Theorem~\ref{EhrTh} will involve transforming $\P(m,n)$ to a unimodularly equivalent generalized permutohedron $\PP(m,n)$ in $\R^{m+1}$ (see~\eqref{Plift}),
obtaining a decomposition of $\PP(m,n)|_{\,n\ge m-1}$ as a Minkowski sum of dilated coordinate simplices (see~\eqref{Mink}),
applying a result of Postnikov~\cite[Thm.~11.3]{Pos09} for the number of integer points in certain generalized permutohedra,
observing that this gives an expression for $\ehr_{\P(m,n)}(t)|_{\,n\ge m-1}$ as an edge-weighted sum
over graphs (with loops and multiple edges permitted) on~$m$ labelled vertices in which each connected component contains at most one cycle (see~\eqref{EhrGra}), and
finally applying standard techniques for the enumeration of such graphs to obtain generating function expressions (see~\eqref{EhrT} and~\eqref{EhrC}) for $\ehr_{\P(m,n)}(t)|_{\,n\ge m-1}$,
from which~\eqref{EhrEq} follows easily by evaluating power series coefficients.

An outline of the remaining sections of this paper is as follows.  In Section~\ref{prop}, a review is provided of key properties of~$\P(m,n)$, including
characterizations of its vertices, facets and other faces.
In Section~\ref{conn}, a review is provided of connections between~$\P(m,n)$ and certain other polytopes.
In Section~\ref{proof}, the proof of Theorem~\ref{EhrTh} is given, where the first step of the proof uses the connection between $\P(m,n)$ and generalized permutohedra described in~(i) of Section~\ref{conn}.
In Section~\ref{disc}, various implications of Theorem~\ref{EhrTh} are discussed, including some related to the connections between $\P(m,n)$ and polytopes described in 
Section~\ref{conn}.

\section{Properties of $\P(m,n)$}\label{prop}
Some known properties of~$\P(m,n)$ are as follows.
\begin{list}{(\roman{c-prop})}{\usecounter{c-prop}\setlength{\labelwidth}{9.1mm}\setlength{\leftmargin}{9.1mm}\setlength{\labelsep}{1.5mm}\setlength{\topsep}{0.9mm}\setlength{\itemsep}{1mm}}
\item It is noted in~\cite[Rem.~5.5]{HeuStr22} that $\P(m,n)$ has dimension $m$.
\item It is shown in~\cite[Prop.~3.5]{BehCasChaDiaEscHarInk22} that $\P(m,n)$ is a simple polytope.
\item It is shown in~\cite[Prop.~5.7]{HeuStr22} that $\P(m,n)$ has $\sum_{i=0}^{\min(m,n)}\frac{m!}{(m-i)!}$ vertices, given explicitly as
the vectors in $\R^m$ with entries of zero in any $m-i$ positions, and with the
other $i$ entries being $n,n-1,\ldots,n-i+1$ in any order, for $i=0,\ldots,\min(m,n)$.
\item It is shown in~\cite[Thms.~5.10 \& 5.11]{HeuStr22} that $\P(m,n)$ can be expressed a bounded intersection of affine halfspaces as
\begin{equation}\label{facets}\P(m,n)=\left\{x\in\R^m\left|\,\;\begin{array}{@{}l@{}}x_i\ge0,\text{ for }i=1,\ldots,m,\\[1.5mm]
\displaystyle\sum_{i\in S}x_i\le\sum_{i=n-|S|+1}^n\!i,\text{ for all }\emptyset\subsetneq S\subsetneq\{1,\ldots,m\}\\[-3mm]
\hspace{34.5mm}\text{ with }|S|\le\min(m,n)-1,\\[2.5mm]
\displaystyle\sum_{i=1}^m x_i\le\sum_{i=\max(1,n-m+1)}^n\!i\end{array}\right.\right\}.\end{equation}
Furthermore, if any single inequality in the set on the RHS of~\eqref{facets} is changed to an equality, then the set is a facet of~$\P(m,n)$,
and all facets arise uniquely in this way,
from which it follows that~$\P(m,n)$ has $m+\sum_{i=0}^{\min(m,n)-1}\binom{m}{i}$ facets.
\item It is noted in~\cite[Eq.~(4.12)]{BehCasChaDiaEscHarInk22} that for $n\ge m-1$, $\P(m,n)$ has the Minkowski sum decomposition
\begin{equation}\label{Mink1}\P(m,n)\big|_{n\ge m-1}=(n-m+1)\sum_{i=1}^m\conv\{0,e_i\}+\sum_{1\le i<j\le m}\conv\{0,e_i,e_j\},\end{equation}
where $\conv$ denotes the convex hull, 
and $e_i$ denotes the $i$th standard unit vector in~$\R^m$. 
Note that the term $\sum_{i=1}^m\conv\{0,e_i\}$ in~\eqref{Mink1} is simply the unit hypercube $[0,1]^m$.
\item It is shown in~\cite[Thm.~3.6]{BehCasChaDiaEscHarInk22} and~\cite[Thm.~7.5]{BlaSan22} that there is a simple
bijection between the $i$-dimensional faces of~$\P(m,n)$ and certain chains of subsets of $\{1,\dots,m\}$ with~$i$ so-called missing
ranks. 
The $m=n$ case of this result was obtained in~\cite[Thm.~5.24]{HeuStr22}, and the general case was conjectured in~\cite[Conj.~5.25]{HeuStr22}.
\item It is shown in~\cite[Thm.~3.19]{BehCasChaDiaEscHarInk22} that the $f$-polynomial $f_{\P(m,n)}(t)$ of $\P(m,n)$ (i.e., the polynomial in~$t$ whose
coefficient of $t^i$ is the number of $i$-dimensional faces of $\P(m,n)$) is given by
\begin{equation*}f_{\P(m,n)}(t)=1+\sum_{i=0}^{n-1}\binom{m}{i}\,A_i(t+1)\:\sum_{j=1}^{m-i}(t+1)^j,\end{equation*}
where $A_i(t)$ is the Eulerian polynomial (defined for $i\ge1$ as the polynomial in~$t$ whose coefficient of $t^j$ is the number of permutations of~$\{1,\ldots,i\}$ with exactly $j$ descents,
and for $i=0$ as $A_0(t)=1$).
\item 
It is shown in~\cite[Cor.~3.11 \& Cor.~3.20]{BehCasChaDiaEscHarInk22} that, for fixed $m$, all $\P(m,n)$ with $n\ge m$ are combinatorially equivalent,
with $f$-polynomial given by
\begin{equation*}f_{\P(m,n)}(t)\big|_{n\ge m-1}=1+(t+1)\,\sum_{i=1}^m\binom{m}{i}A_i(t+1).\end{equation*}
\item It is shown in~\cite[Thm.~4.5]{BehCasChaDiaEscHarInk22} and~\cite[Cor.~3.28]{HanLenVin23} that for $n\ge m-1$,
the volume $\vol\P(m,n)$ of~$\P(m,n)$ is explicitly
\begin{equation}\label{volP}\vol\P(m,n)\big|_{n\ge m-1}=-\frac{1}{2^m}\,\sum_{i=0}^m\,\binom{m}{i}(2i-3)!!\,(2n+1)^{m-i}.\end{equation}
Related results are obtained for $n=m-1$ in~\cite[Sec.~4]{AmaWan22} and~\cite[Part~(d)]{Sto22}.
\item Explicit expressions are obtained for the Ehrhart polynomial of~$\P(m,n)$ with arbitrary~$m$ and fixed $n\le 3$
in~\cite[p.~28, Thm.~5.11~\& Thm.~5.12]{BehCasChaDiaEscHarInk22}, and for the volume of~$\P(m,4)$ with arbitrary~$m$ in~\cite[Thm.~5.13]{BehCasChaDiaEscHarInk22}.
Explicit expressions are obtained for the Ehrhart polynomial of~$\P(m,n)$ with fixed $m\le4$ and arbitrary $n\ge m-1$ in~\cite[Eq.~(6.1), Eq.~(6.2), Thm.~6.1 \& Thm.~6.2]{BehCasChaDiaEscHarInk22},
and these expressions can now be obtained alternatively using~\eqref{EhrEq} with $m\le 4$.
\end{list}

\section{Connections between $\P(m,n)$ and other polytopes}\label{conn}
Some connections between partial permutohedra and certain other polytopes are as follows.
\begin{list}{(\roman{c-conn})}{\usecounter{c-conn}\setlength{\labelwidth}{9.1mm}\setlength{\leftmargin}{9.1mm}\setlength{\labelsep}{1.5mm}\setlength{\topsep}{0.9mm}\setlength{\itemsep}{1mm}}
\item\emph{Generalized permutohedra}.
Let $\PP(m,n)$ be the polytope which is obtained by lifting $\P(m,n)$ from $\R^m$ into the hyperplane $\{x\in\R^{m+1}\mid\sum_{i=1}^{m+1}x_i=\sum_{i=\max(1,n-m+1)}^ni\}$ in $\R^{m+1}$
according to
\begin{equation}\label{Plift}\textstyle\PP(m,n)=\Bigl\{\bigl(x_1,\ldots,x_m,\sum_{i=\max(1,n-m+1)}^ni-\sum_{i=1}^mx_i\bigr)\Bigm|x\in\P(m,n)\Bigr\}.\end{equation}
It follows that $\PP(m,n)$ is unimodularly equivalent to~$\P(m,n)$, and hence that it has the same Ehrhart polynomial as~$\P(m,n)$.

It is noted in~\cite[Sec.~4.3]{BehCasChaDiaEscHarInk22} that $\PP(m,n)$ is a case of a generalized permutohedron, as defined by Postnikov~\cite[Def.~6.1]{Pos09}.
Furthermore, as noted in~\cite[Eq.~(4.14)]{BehCasChaDiaEscHarInk22}, it follows from~\eqref{Mink1} and~\eqref{Plift} that for $n\ge m-1$, $\PP(m,n)$ has the Minkowski sum decomposition
\begin{equation}\label{Mink}\PP(m,n)\big|_{n\ge m-1}=(n-m+1)\sum_{i=1}^m\conv\{\e_i,\e_{m+1}\}+\sum_{1\le i<j\le m}\conv\{\e_i,\e_j,\e_{m+1}\},\end{equation}
where $\e_i$ denotes the $i$th standard unit vector in~$\R^{m+1}$.  It can be seen from~\eqref{Mink} that for $n\ge m-1$,~$\PP(m,n)$ is
a type-$\mathcal{Y}$ generalized permutohedron,
i.e., it has the form $P^y_n(\{y_I\})$, as defined in~\cite[p.~1042]{Pos09}.
\item\emph{Antiblocking permutohedra}.  For $z_1,\ldots,z_m\in\R$, the standard permutohedron $\Pi(z_1,\ldots,z_m)$ is defined as the convex hull of all
vectors whose entries are permutations of $z_1,\ldots,z_m$, i.e.,
\begin{equation*}\Pi(z_1,\ldots,z_m)=\conv\{(z_{\sigma(1)},\ldots,z_{\sigma(m)})\mid\sigma\in\mathfrak{S}_m\},\end{equation*}
 and for nonnegative $z_1,\ldots,z_m$, the antiblocking permutohedron $\widehat{\Pi}(z_1,\ldots,z_m)$ is defined in~\cite[Def.~5.3]{BehCasChaDiaEscHarInk22} as
\begin{multline*}\hspace{10mm}\widehat{\Pi}(z_1,\ldots,z_m)\\
=\bigl\{x\in\R^m\bigm|\text{there exists }y\in\Pi(z_1,\ldots,z_m)\text{ with }0\le x_i\le y_i\text{ for }i=1,\ldots,m\bigr\}.\end{multline*}
As noted in~\cite[Eqs.~5.4 \&~5.5]{BehCasChaDiaEscHarInk22}, $\widehat{\Pi}(z_1,\ldots,z_m)$ can also be expressed as
\begin{align*}\widehat{\Pi}(z_1,\ldots,z_m)&=\bigl\{(x_1y_1,\ldots,x_my_m)\bigm|x\in[0,1]^m,\ y\in\Pi(z_1,\ldots,z_m)\bigr\},\\
&=\conv\bigl\{(a_1z_{\sigma(1)},\ldots,a_mz_{\sigma(m)})\bigm|a\in\{0,1\}^m,\ \sigma\in\mathfrak{S}_m\bigr\}.\end{align*}
It is shown in~\cite[Cor.~5.8]{BehCasChaDiaEscHarInk22} that $\P(m,n)$ is an antiblocking permutohedron, given by
\begin{equation}\label{Pab}\P(m,n)=\begin{cases}
\widehat{\Pi}(n,n-1,n-2,\dots,1,\underbrace{0,\ldots,0\!}_{m-n}\,),&\text{ for }n\le m-2,\\[4mm]
\widehat{\Pi}(n,n-1,n-2,\dots,n-m+1),&\text{ for }n\ge m-1.\end{cases}\end{equation}
Note that the permutohedron $\Pi(n,n-1,n-2,\dots,1,\underbrace{0,\ldots,0\!}_{m-n}\,)$ for $n\le m-2$,
or $\Pi(n,n-1,n-2,\dots,n-m+1)$ for $n\ge m-1$, is a facet of $\P(m,n)$, as given
by the set on the RHS of~\eqref{facets} with its last inequality changed to an equality.
\item\emph{Parking function polytopes}. A parking function of length~$m$ is an $m$-vector of positive integers whose nondecreasing rearrangement $r_1\le r_2\le\ldots\le r_m$ satisfies $r_i\le i$ for $i=1,\ldots,m$,
and the parking function polytope~$P_m$, as defined in~\cite{Sta20}, is the convex hull of all parking functions of length~$m$.
The vertices and facets of~$P_m$ are characterized in~\cite[Sec.~1]{AmaWan22} and~\cite[Parts (a) \&~(b)]{Sto22},
and by comparing these with the vertices or facets of~$\P(m,m-1)$, it follows that~$P_m$ with $m\ge2$ is simply~$\P(m,m-1)$ translated by $(1,\ldots,1)$, i.e.,
\begin{equation}\label{park}\P(m,m-1)=\bigl\{(x_1-1,\ldots,x_m-1)\mid x\in P_m\bigr\}.\end{equation}
A generalization of~$P_m$ is defined in~\cite{HanLenVin23} as follows.
For positive integers~$a$,~$b$ and~$m$, an $(a,b)$-parking function
of length~$m$ is an $m$-vector of positive integers whose nondecreasing rearrangement $r_1\le r_2\le\ldots\le r_m$ satisfies $r_i\le a+(i-1)b$ for $i=1,\ldots,m$, and the
$(a,b)$-parking function polytope~$\mathcal{X}_m(a,b)$ is defined as the convex hull of all $(a,b)$-parking functions of length~$m$.
It follows that $P_m=\mathcal{X}_m(1,1)$, and it is shown in~\cite[Prop.~3.16]{HanLenVin23} that for $n\ge m-1$, $\mathcal{X}_m(n-m+2,1)$ is~$\P(m,n)$ translated by $(1,\ldots,1)$, i.e.,
\begin{equation}\label{genpark}\P(m,n)\big|_{n\ge m-1}=\bigl\{(x_1-1,\ldots,x_m-1)\mid x\in\mathcal{X}_m(n-m+2,1)\bigr\}.\end{equation}
\item\emph{Win vector polytopes}.
For a graph~$G$ with vertices $1,\ldots,m$, a partial orientation~$\mathcal{O}$ of~$G$ is an assignment of a direction to some edges of~$G$,
and the win vector of~$\mathcal{O}$ is the indegree sequence of~$\mathcal{O}$, i.e., for $i=1\,\ldots,m$, the $i$th entry of the win vector is the number of edges incident to~$i$ which are
directed towards~$i$ by~$\mathcal{O}$. The win vector polytope~$W(G)$ of~$G$,
as defined in~\cite{BarMouWel97}, is then the convex hull 
of the win vectors of all partial orientations of~$G$.
For the complete graph~$K_m$ with $m\ge2$,
\begin{equation}\label{win}\P(m,m-1)=W(K_m),\end{equation}
where this can be shown using a characterization of the vertices of win vector polytopes given in~\cite[Prop.~3.4]{BarMouWel97}.
(although note that the description of the win vectors of~$K_m$ given in~\cite[Ex.~3]{BarMouWel97} seems to contain errors).
\item\emph{Stellohedra}. For a graph~$G$ with vertex set~$1,\ldots,N$, the graph associahedron $\textrm{Assoc}(G)$ of~$G$ is
defined as $\textrm{Assoc}(G)=\sum_S\conv\{e_i\mid i\in S\}$,
where this is a Minkowski sum over all nonempty and nonsingleton subsets~$S$ of~$\{1,\ldots,N\}$ such that the subgraph of~$G$ induced by $S$ is connected, and $e_i$ denotes the $i$th standard unit vector in $\R^N$.
The $m$-stellohedron, as defined in~\cite[Sec.~10.4]{PosReiWil08}, is then the graph associahedron
of the star graph $K_{1,m}$, consisting of a central vertex $m+1$ connected to~$m$ vertices~$1,\ldots,m$, which gives
\begin{equation*}\textrm{Assoc}(K_{1,m})=\sum_{\emptyset\subsetneq S\subseteq\{1,\ldots,m\}}\conv\bigl(\{\e_i\mid i\in S\}\cup\{\e_{m+1}\}\bigr),\end{equation*}
where $\e_i$ denotes the $i$th standard unit vector in $\R^{m+1}$.
It is noted in~\cite[Thm.~5.17]{HeuStr21} that~$\P(m,m)$
is combinatorially equivalent to the $m$-stellohedron.  Hence, due to the combinatorial equivalence discussed in~(viii) in Section~\ref{prop}, all
$\P(m,n)$ with $n\ge m$ are combinatorially equivalent to the $m$-stellohedron, i.e.,
\begin{equation}\label{stell}\P(m,n)\big|_{n\ge m}\cong\textrm{Assoc}(K_{1,m}),\end{equation}
where $\cong$ denotes combinatorial equivalence.  For further information, see for example~\cite[Rem.~3.12]{BehCasChaDiaEscHarInk22}.
\item\emph{Monotone path polytopes of polymatroids}.
As shown in~\cite[Thm.~7.3]{BlaSan22}, $\P(m,n)$ is also connected to certain so-called monotone path polytopes of polymatroids.
As noted in~\cite[Cor.7.4]{BlaSan22}, this provides a generalization to all~$m$ and~$n$ of the combinatorial equivalence of~\eqref{stell}, specifically
\begin{equation*}\P(m,n)\cong
\sum_{\substack{S\subseteq[m]\\|S|\ge\max(1,m-n+1)}}\!\!\!\conv\bigl(\{\e_i\mid i\in S\}\cup\{\e_{m+1}\}\bigr).\end{equation*}
For further information, see also~\cite[Rem.~3.13]{BehCasChaDiaEscHarInk22}.
\end{list}

\section{Proof of Theorem~\ref{EhrTh}}\label{proof}
A proof of Theorem~\ref{EhrTh} is provided in this section.
The proof is subdivided into three parts, and an outline of these is as follows.
\begin{list}{$\bullet$}{\setlength{\labelwidth}{3mm}\setlength{\leftmargin}{5mm}\setlength{\labelsep}{2mm}\setlength{\topsep}{0.9mm}\setlength{\itemsep}{1mm}}
\item In Subsection~\ref{PosApp}, a certain finite set~$\A(m)$ of nonnegative integer sequences is defined first.
A theorem of Postnikov~\cite[Thm.~11.3]{Pos09}
for the number of integer points in any type-$\mathcal{Y}$ generalized permutohedron is then applied to
the Minkowski sum decomposition~\eqref{Mink} of the lifted version~$\PP(m,n)$ of $\P(m,n)$ with $n\ge m-1$,
and this leads to an expression (see~\eqref{EhrDrac}) for $\ehr_{\P(m,n)}(t)|_{n\ge m-1}$ as a
sum over~$\A(m)$.
\item In Subsection~\ref{graphs}, properties of connected graphs with at most one cycle are considered first, and a set~$\G(m)$ is then defined as the set of
all graphs with vertices labelled $1,\ldots,m$, such that loops and multiple edges are permitted, and
each connected component contains at most one cycle.
The set~$\A(m)$ is then found to be in simple bijection with~$\G(m)$, and applying this
bijection to the previous expression~\eqref{EhrDrac} for $\ehr_{\P(m,n)}(t)|_{n\ge m-1}$ gives a further expression (see~\eqref{EhrGra})
for~$\ehr_{\P(m,n)}(t)|_{n\ge m-1}$ as a certain edge-weighted sum over the graphs of~$\G(m)$.
\item In Subsection~\ref{graphenum}, an exponential generating function for $\ehr_{\P(m,n)}(t)|_{n\ge m-1}$ is introduced (see~\eqref{EhrFDef}),
and using the expression~\eqref{EhrGra} for $\ehr_{\P(m,n)}(t)|_{n\ge m-1}$ then enables this generating function to be expressed in terms
of exponential generating functions for certain graphs (see~\eqref{fptz}-\eqref{fptz2}).  The required expression~\eqref{EhrEq} for $\ehr_{\P(m,n)}(t)|_{n\ge m-1}$
follows by applying standard techniques for obtaining explicit expressions for such graph-related generating functions.
\end{list}

\subsection{Application of Postnikov's theorem}\label{PosApp}
Let $\E(m)$ be the set of all subsets of $\{1,\ldots,m\}$ of size~1 or~2, i.e.,
\begin{equation*}\E(m)=\bigl\{\{i\}\bigm|i=1,\ldots,m\bigr\}\cup\bigl\{\{i,j\}\bigm|1\le i<j\le m\bigr\},\end{equation*}
and let $\A(m)$ be the set of all nonnegative integer sequences $a=(a_S)_{S\in\E(m)}$
which satisfy
\begin{equation}\label{Dm}\sum_{S\in E}a_S\le\Bigl|\bigcup_{S\in E}S\Bigr|\text{ \ for all }\emptyset\subsetneq E\subseteq\E(m).\end{equation}
As examples of the condition in~\eqref{Dm}, by taking $E=\bigl\{\{i\}\bigr\}$ for $i=1,\ldots,m$, $E=\bigl\{\{i,j\}\bigr\}$ for $1\le i<j\le m$, and $E=\E(m)$ in~\eqref{Dm}, it follows that any $a\in\A(m)$ satisfies
\begin{equation*}a_{\{i\}}\in\{0,1\}\text{ for }i=1,\ldots,m,\ \ a_{\{i,j\}}\in\{0,1,2\}\text{ for }1\le i<j\le m,\text{ \ and }\sum_{S\in\E(m)}\!a_S\le m.\end{equation*}
As an example of a set $\A(m)$, if $m=2$, then
\begin{equation}\label{A2}\A(2)=\bigl\{(0,0,0),\,(0,0,1),\,(0,0,2),\,(1,0,0),\,(0,1,0),\,(1,1,0),\,(1,0,1),\,(0,1,1)\bigr\},\end{equation}
where each $a\in\A(2)$ is written as $(a_{\{1\}},a_{\{2\}},a_{\{1,2\}})$.

The condition~\eqref{Dm} can also be interpreted in terms of systems of distinct representatives, as follows.
For a nonnegative integer sequence $a=(a_S)_{S\in\E(m)}$,  
let $\Phi(a)$ be the family (or multiset) of $\sum_{S\in\E(m)}a_S$ elements of $\E(m)$ in which~$S$ appears exactly $a_S$ times
for each $S\in\E(m)$.  It then follows, by applying Hall's marriage theorem to~$\Phi(a)$,
that
\begin{equation}\label{SDR}\Phi(a)\text{ has a system of distinct representatives if and only if }a\text{ satisfies }\eqref{Dm}.\end{equation}
Note that if $a_S=0$ for each $S\in\E(m)$, then $\Phi(a)$ is the empty family, which is considered to have a system of distinct representatives.

It will now be shown, using Postnikov's theorem~\cite[Thm.~11.3]{Pos09}, that
\begin{equation}\label{EhrDrac}
\ehr_{\P(m,n)}(t)\big|_{n\ge m-1}=\sum_{a\in\A(m)}\:\prod_{i=1}^m\binom{(n-m+1)t+a_{\{i\}}-1}{a_{\{i\}}}
\:\prod_{1\le i<j\le m}\binom{t+a_{\{i,j\}}-1}{a_{\{i,j\}}},\end{equation}
where this expression is also given in~\cite[Eq.~(6.5)]{BehCasChaDiaEscHarInk22}, and the details of its derivation
(which are omitted in~\cite{BehCasChaDiaEscHarInk22}) are as follows.
For $n\ge m-1$, let $\mathcal{Q}=\conv\{\e_1,\ldots,\e_{m+1}\}+t\,\PP(m,n)$, using the same notation
as in~\eqref{Mink}.   Then the trimmed version (as defined in~\cite[Def.~11.2]{Pos09}) of $\mathcal{Q}$
is $\mathcal{Q}^{-}=t\,\PP(m,n)$.  Applying~\cite[Thm.~11.3]{Pos09} to~$\mathcal{Q}^{-}$, decomposing $\PP(m,n)$ according to~\eqref{Mink},
using the fact that $\ehr_{\PP(m,n)}(t)$ is the number of integer points in $t\,\PP(m,n)$ for any nonnegative integer $t$,
and recalling (as indicated after~\eqref{Plift}) that $\P(m,n)$ and $\PP(m,n)$ have the same Ehrhart polynomial,
then gives
\begin{equation}\label{EhrDrac1}\ehr_{\P(m,n)}(t)\big|_{n\ge m-1}=\sum_{b\in\A'(m)}\:\prod_{i=1}^m\binom{(n-m+1)t+b_{\{i\}}-1}{b_{\{i\}}}
\:\prod_{1\le i<j\le m}\binom{t+b_{\{i,j\}}-1}{b_{\{i,j\}}},\end{equation}
where a trivial factor $\bigl(\begin{smallmatrix}b_{\{1,\ldots,m+1\}}\\b_{\{1,\ldots,m+1\}}\end{smallmatrix}\bigr)=1$ initially appears
in the sum on the RHS, and
where (due to~\cite[Def.~9.2]{Pos09}) $\A'(m)$ is the set of all nonnegative integer sequences $b=(b_S)_{S\in\E'(m)}$ such that
\begin{equation*}\sum_{S\in\E'(m)}b_S=m,\text{ \ \ and \ \ }
\sum_{S\in E}b_S\le\Bigl|\bigcup_{S\in E}S\Bigr|-1\text{ \ for all }\emptyset\subsetneq E\subsetneq\E'(m),\end{equation*}
with
\begin{equation*}\E'(m)=\bigl\{\{1,\ldots,m+1\}\bigr\}\cup\bigl\{\{i,m+1\}\bigm|i=1,\ldots,m\bigr\}\cup\bigl\{\{i,j,m+1\}\bigm|1\le i<j\le m\bigr\}.\end{equation*}
It can be checked straightforwardly that there is a bijection from $\A'(m)$ to~$\A(m)$, in which $b\in\A'(m)$ is mapped to a sequence in $\A(m)$ by
simply omitting the entry $b_{\{1,\ldots,m+1\}}$.
Hence, the sum over $\A'(m)$ in~\eqref{EhrDrac1} can be replaced by a sum over $\A(m)$, thereby giving~\eqref{EhrDrac}.

As an aside (and as also discussed in \cite[Sec.~6.2)]{BehCasChaDiaEscHarInk22}),
note that certain properties of $\ehr_{\P(m,n)}(t)$ follow immediately from~\eqref{EhrDrac}.
For example, it can be seen that, for fixed $m$ and $t$,
$\ehr_{\P(m,n)}(t)$ is a polynomial in $(n-m+1)t$,
that all coefficients of this polynomial are positive integers if~$t$ is a positive integer,
and that this polynomial has degree~$m$ (since $\A(m)$ contains the sequence~$a$ with entries $a_{\{i\}}=1$ for $i=1,\ldots,m$ and $a_{\{i,j\}}=0$ for $1\le i<j\le m$).
It can also be seen that, for fixed~$m$ and~$n$,
$\ehr_{\P(m,n)}(t)$  is a polynomial in~$t$ of degree~$m$,
where this also follows from general Ehrhart theory,
and that all coefficients of this polynomial are positive, where this also follows from the general property
that the Ehrhart polynomial of any integral type-$\mathcal{Y}$ generalized permutohedron has positive coefficients (see~\cite[Cor.~3.1.5]{Liu19}).
For example, for $m=2$,~\eqref{A2} and~\eqref{EhrDrac} give
\begin{align*}\textstyle\ehr_{\P(2,n)}(t)&\textstyle=1+\binom{t}{1}+\binom{t+1}{2}+\binom{(n-1)t}{1}+\binom{(n-1)t}{1}+\binom{(n-1)t}{1}^2+\binom{(n-1)t}{1}\binom{t}{1}+\binom{(n-1)t}{1}\binom{t}{1}\\
&=(n^2-1/2)\,t^2+(2n-1/2)\,t+1,\end{align*}
which, for this simple case, could also easily be obtained directly (or by using Pick's theorem)
by observing that $\P(2,n)$ is the pentagon~$\conv\{(0,0),\,(n,0),\,(n,n-1),\,(n-1,n),(0,n)\}$.

\subsection{Introduction of graphs}\label{graphs}
In this subsection and Subsection~\ref{graphenum}, undirected graphs with loops and multiple edges permitted are considered.  For such graphs, the following
conventions are used. A cycle of length~1 corresponds to a vertex with a loop attached, a cycle of length~2 corresponds to a pair of distinct vertices connected by two parallel edges,
and a cycle of length~$\ell\ge3$ corresponds to distinct vertices $i_1,\ldots,i_\ell$, such that~$i_k$ and~$i_{k+1}$ are adjacent for $k=1,\ldots,\ell-1$, and~$i_\ell$ and~$i_1$ are adjacent.
An orientation of a graph is an assignment of a direction to all edges, where there is only one choice for the direction of each loop, and two choices for the direction
of each other edge. For such an orientation, the indegree of a vertex $i$ is the number of edges incident to~$i$ which are
directed towards~$i$, where each loop attached to~$i$ contributes~$1$ to this number.

Consider a connected graph $G$, with loops and multiple
edges permitted. 
Then, as will be discussed below, the following are all equivalent:
\begin{list}{(\roman{c-gr})}{\usecounter{c-gr}\setlength{\labelwidth}{13mm}\setlength{\leftmargin}{13mm}\setlength{\labelsep}{1.5mm}\setlength{\topsep}{0.9mm}\setlength{\itemsep}{0.8mm}}
\item $G$ is acyclic, i.e., a tree.
\item The number of edges in $G$ is equal to the number of vertices in $G$ minus~$1$. 
\item There exists an orientation of $G$ in which one vertex has indegree~$0$,
and each other vertex has indegree~$1$.
\end{list}
Similarly, the following are all equivalent:
\begin{list}{(\roman{c-gr}$'$)}{\usecounter{c-gr}\setlength{\labelwidth}{13mm}\setlength{\leftmargin}{13mm}\setlength{\labelsep}{1.5mm}\setlength{\topsep}{0.9mm}\setlength{\itemsep}{0.8mm}}
\item $G$ contains exactly one cycle.
\item The number of edges in $G$ is equal to the number of vertices in $G$. 
\item There exists an orientation of $G$ in which each vertex has indegree~$1$.
\end{list}

The equivalences between~(i) and~(ii), and between~(i$'$) and~(ii$'$), are standard facts in graph theory.  The implications from~(iii) to~(ii),
and from~(iii$'$) to~(ii$'$), follow from the observation that in any oriented graph, the sum of indegrees over all vertices equals the number of edges.
The implication from~(i) to~(iii) can be confirmed by choosing any vertex~$i$ of a tree~$G$, and for an edge which connects a vertex~$j$ of distance~$d$ from~$i$
to a vertex~$j'$ of distance $d+1$ from~$i$, assigning the direction from~$j$ to~$j'$.
The implication from~(i$'$) to~(iii$'$) can be confirmed as follows.  First observe that if~$G$ contains exactly one cycle, with vertices $i_1,\ldots,i_\ell$,
then~$G$ consists of that cycle together with~$\ell$ mutually disjoint trees $\tau_1,\ldots,\tau_\ell$,
such that $i_k$ is a vertex of~$\tau_k$, for $k=1,\ldots,\ell$.  The required orientation of such a graph can then be obtained by directing the edges of the cycle to form a directed path around the cycle
(which can be done in two ways for $\ell\ge2$), and for an edge which connects a vertex~$j$ in~$\tau_k$ of distance~$d$ from~$i_k$
to a vertex~$j'$ of distance $d+1$ from~$i_k$ (where these distances are taken within~$\tau_k$), assigning the direction from~$j$ to~$j'$, for $k=1,\ldots,\ell$.


Now let~$\G(m)$ be the set of graphs with vertices labelled $1,\ldots,m$ (and with loops and multiple edges permitted), such that
each connected component contains at most one cycle.

For example, for $m=2$,
\psset{unit=1.3mm}\begin{equation}\label{G2}\G(2)=\bigl\{\pspicture(-4,0)(4.6,0)\multirput(-2,0)(4,0){2}{$\scriptstyle\bullet$}\rput[r](-2.7,0){$\scriptstyle1$}\rput[l](2.8,0){$\scriptstyle2$}\endpspicture,
\pspicture(-5.4,0)(4.6,0)\psline[linewidth=0.5pt](-2,0)(2,0)\multirput(-2,0)(4,0){2}{$\scriptstyle\bullet$}\rput[r](-2.7,0){$\scriptstyle1$}\rput[l](2.8,0){$\scriptstyle2$}\endpspicture,
\pspicture(-5.4,0)(4.6,0)\psellipse[linewidth=0.5pt](0,0)(2,1.5)\multirput(-2,0)(4,0){2}{$\scriptstyle\bullet$}\rput[r](-2.7,0){$\scriptstyle1$}\rput[l](2.8,0){$\scriptstyle2$}\endpspicture,
\pspicture(-5.4,0)(4.6,0)\psellipse[linewidth=0.5pt](-2,1.2)(0.7,1.2)\multirput(-2,0)(4,0){2}{$\scriptstyle\bullet$}\rput[r](-2.7,0){$\scriptstyle1$}\rput[l](2.8,0){$\scriptstyle2$}\endpspicture,
\pspicture(-5.4,0)(4.6,0)\psellipse[linewidth=0.5pt](2,1.2)(0.7,1.2)\multirput(-2,0)(4,0){2}{$\scriptstyle\bullet$}\rput[r](-2.7,0){$\scriptstyle1$}\rput[l](2.8,0){$\scriptstyle2$}\endpspicture,
\pspicture(-5.4,0)(4.6,0)\psellipse[linewidth=0.5pt](-2,1.2)(0.7,1.2)\psellipse[linewidth=0.5pt](2,1.2)(0.7,1.2)\multirput(-2,0)(4,0){2}{$\scriptstyle\bullet$}\rput[r](-2.7,0){$\scriptstyle1$}\rput[l](2.8,0){$\scriptstyle2$}\endpspicture,
\pspicture(-5.4,0)(4.6,0)\psline[linewidth=0.5pt](-2,0)(2,0)\psellipse[linewidth=0.5pt](-2,1.2)(0.7,1.2)\multirput(-2,0)(4,0){2}{$\scriptstyle\bullet$}\rput[r](-2.7,0){$\scriptstyle1$}\rput[l](2.8,0){$\scriptstyle2$}\endpspicture,
\pspicture(-5.4,0)(4,0)\psline[linewidth=0.5pt](-2,0)(2,0)\psellipse[linewidth=0.5pt](2,1.2)(0.7,1.2)\multirput(-2,0)(4,0){2}{$\scriptstyle\bullet$}\rput[r](-2.7,0){$\scriptstyle1$}\rput[l](2.8,0){$\scriptstyle2$}\endpspicture\bigr\}.\end{equation}

It will now be shown that there exists a bijection $\Gamma$ from~$\A(m)$ (as defined in Subsection~\ref{PosApp}) to~$\G(m)$,
where for $a\in\A(m)$, $\Gamma(a)$ is the graph with vertices labelled $1,\ldots,m$, and with
$a_{\{i\}}$ loops attached to vertex $i$ for $i=1,\ldots,m$, and $a_{\{i,j\}}$ edges connecting vertices~$i$ and~$j$ for $1\le i<j\le m$.

Note that the elements of $\G(2)$ in~\eqref{G2} are listed in the order which corresponds, using~$\Gamma$, to the order in which the elements of~$\A(2)$ are listed in~\eqref{A2}.

To confirm the well-definedness of~$\Gamma$, consider $a\in\A(m)$, and associate each edge of~$\Gamma(a)$ with its set of endpoints
(so that loops are associated with singletons, and other edges are associated with sets of size~2). Then~$\Phi(a)$ 
is the family of edges of~$\Gamma(a)$, and due to~\eqref{SDR}, $\Phi(a)$ has a system~$s$ of distinct representatives.
Now form an orientation of~$\Gamma(a)$ by
directing each edge~$e$ towards the endpoint of~$e$ which represents~$e$ in~$s$.  It follows that
for this orientation, each vertex has indegree at most~$1$.  Hence, due to the implications above from~(iii) to~(i),
and from~(iii$'$) to~(i$'$), each connected component of $\Gamma(a)$ contains at most one cycle, so $\Gamma(a)\in\G(m)$, as required.

It is clear from the definition of~$\Gamma$ that any~$a\in\A(m)$ can be uniquely recovered from $\Gamma(a)$, and hence $\Gamma$ is injective.

To confirm the surjectivity of~$\Gamma$, consider $G\in\G(m)$, associate each edge of~$G$ with its set of endpoints, and
form a nonnegative integer sequence $a=(a_S)_{S\in\E(m)}$ by setting~$a_S$ to be the number of edges~$S$ in~$G$, for each $S\in\E(m)$.  Hence,~$\Phi(a)$ is the family of edges of~$G$.
Due to the implications above from~(i) to~(iii),
and from~(i$'$) to~(iii$'$), there exists an orientation~$\mathcal{O}$ of~$G$ in which each vertex has indegree at most~$1$.
A system of distinct representatives of $\Phi(a)$ 
can then be obtained by representing each edge~$e$ of~$G$ by the endpoint of~$e$ towards which~$e$ is directed in~$\mathcal{O}$.
It follows, using~\eqref{SDR}, that $a\in\A(m)$ and $G=\Gamma(a)$, as required.


For $G\in\G(m)$, define an edge pair in $G$ to be a pair of parallel edges connecting the same two vertices, and a single edge in~$G$ to be an edge which is neither a loop nor an edge within an edge pair,
and let $\Nloops{G}$, $\Nsing{G}$ and $\Ndoub{G}$ denote the number of loops, single edges and edge pairs, respectively, in $G$.

Applying the bijection $\Gamma$ to $\A(m)$ in~\eqref{EhrDrac} 
now gives
\begin{multline}\label{EhrGra}
\ehr_{\P(m,n)}(t)\big|_{n\ge m-1}\\
=\sum_{G\in\G(m)}\:\bigl((n-m+1)t\bigr)^{\Nloops{\,G}}\:t^{\Nsing{\,G}}\:\bigl(t(t+1)/2\bigr)^{\Ndoub{\,G}}.\end{multline}

\subsection{Graph enumeration}\label{graphenum}
In this subsection, 
$\P(m,n)$ with $n\ge m-1$ is initially considered in the form $\P(m,m+p-1)$,
where~$p$ is a fixed nonnegative integer.

Consider the exponential generating function $F_p(t,z)$ for $\ehr_{\P(m,m+p-1)}(t)$, as defined by
\begin{equation}\label{EhrFDef}F_p(t,z)=1+\sum_{m=1}^\infty\,\ehr_{\P(m,m+p-1)}(t)\,\frac{z^m}{m!}.\end{equation} 
Substituting~\eqref{EhrGra} into~\eqref{EhrFDef} gives
\begin{equation*}F_p(t,z)
=1+\sum_{m=1}^\infty\,\sum_{G\in\G(m)}
(p\,t)^{\Nloops{\,G}}\:t^{\Nsing{\,G}}\:\bigl(t(t+1)/2\bigr)^{\Ndoub{\,G}}\,\frac{z^m}{m!}.\end{equation*}
Now let $\Gc(m)$ denote the set of connected graphs in $\G(m)$, and define
\begin{equation}\label{fptz}f_p(t,z)=\sum_{m=1}^\infty\,\sum_{G\in\Gc(m)}(p\,t)^{\Nloops{\,G}}\:t^{\Nsing{\,G}}\:\bigl(t(t+1)/2\bigr)^{\Ndoub{\,G}}\,\frac{z^m}{m!}.\end{equation}
It follows from the exponential formula for exponential generating functions (see for example~\cite[Secs.~5.1--5.3]{Sta99}) that
\begin{equation}\label{Fptz}F_p(t,z)=e^{f_p(t,z)}.\end{equation}

Since $\G(m)$ consists of graphs in which each connected component contains at most one cycle, a graph in~$\Gc(m)$ is one of the following.
\begin{list}{$\bullet$}{\setlength{\labelwidth}{5mm}\setlength{\leftmargin}{7mm}\setlength{\labelsep}{2mm}\setlength{\topsep}{0.9mm}\setlength{\itemsep}{0.8mm}}
\item A tree. 
\item A graph obtained by attaching a loop to one vertex of a tree. Such a graph will be referred to as a looped tree.
Note that a looped tree can be associated with a rooted tree by regarding the vertex with the loop attached as the root, and then deleting the loop.
\item A graph obtained by converting one edge of a tree to an edge pair.  Such a graph will be referred to as an enhanced tree.
\item A connected graph without any loops or multiple edges, and with exactly one cycle.  Such a graph will be referred to as a quasitree.
\end{list}
Denote the associated subsets of $\Gc(m)$ as $\Gtree(m)$, $\Gltree(m)$, $\Getree(m)$ and $\Gqtree(m)$, respectively.

The implications from~(i) to~(ii), and from~(i$'$) to~(ii$'$) in Subsection~\ref{graphs}, give the following data for $G\in\Gc(m)$:
\begin{equation*}\begin{array}{c|c|c|c}&\,\Nloops{G}&\,\Nsing{G}&\,\Ndoub{G}\\[0.8mm]\hline
\rule{0mm}{5.5mm}G\in\Gtree(m)&0&m-1&0\\[1mm]
G\in\Gltree(m)&1&m-1&0\\[1mm]
G\in\Getree(m)&0&m-2&1\\[1mm]
G\in\Gqtree(m)&0&m&0\end{array}\end{equation*}

Expressing $\Gc(m)$ in~\eqref{fptz} as the disjoint union of $\Gtree(m)$, $\Gltree(m)$, $\Getree(m)$ and $\Gqtree(m)$ 
then gives
\begin{multline}\label{fptz2}
f_p(t,z)=\frac{1}{t}\,\sum_{m=1}^\infty|\Gtree(m)|\,\frac{(tz)^m}{m!}+p\,\sum_{m=1}^\infty|\Gltree(m)|\,\frac{(tz)^m}{m!}\\
+\frac{t+1}{2t}\,\sum_{m=1}^\infty|\Getree(m)|\,\frac{(tz)^m}{m!}+\sum_{m=1}^\infty|\Gqtree(m)|\,\frac{(tz)^m}{m!}.
\end{multline}

Now consider the function
\begin{equation}\label{Tdef}T(z)=\sum_{m=1}^\infty m^{m-1}\frac{z^m}{m!}.\end{equation}
For information on $T(z)$, see for example~\cite[pp.~331, 332 and~338]{CorGonHarJefKnu96} or~\cite[pp.~23--28 and~43]{Sta99}.
Note that $T(z)$ is related to the Lambert~$W$ function $W(z)$ by $T(z)=-W(-z)$.

Cayley's formula for the number of trees on $m$ labelled vertices states that $|\Gtree(m)|=m^{m-2}$.  Therefore,
$|\Gltree(m)|=m^{m-1}$, since every looped tree in $\Gltree(m)$ can be uniquely obtained by choosing any tree in $\Gtree(m)$,
and choosing any one of the~$m$ vertices of the tree at which to attach a loop.  Similarly,
$|\Getree(m)|=(m-1)m^{m-2}$, since every enhanced tree in $\Getree(m)$ can be uniquely obtained by choosing any tree in $\Gtree(m)$,
and choosing any one of the~$m-1$ edges of the tree to be converted to an edge pair.

It follows immediately using~\eqref{Tdef} that
\begin{equation}\label{Fltree}\sum_{m=1}^\infty|\Gltree(m)|\,\frac{z^m}{m!}=T(z),\end{equation}
and it can be shown that
\begin{align}\label{Fetree}\sum_{m=1}^\infty|\Getree(m)|\,\frac{z^m}{m!}&=\frac{T(z)^2}{2},\\
\label{Ftree}\sum_{m=1}^\infty|\Gtree(m)|\,\frac{z^m}{m!}&=T(z)-\frac{T(z)^2}{2},\end{align}
where~\eqref{Fetree} can be obtained by observing that every enhanced tree corresponds uniquely to a set $\{\tau_1,\tau_2\}$ of disjoint
rooted trees~$\tau_1$ and~$\tau_2$ (since an enhanced tree can be viewed as~$\tau_1$ and~$\tau_2$, together with an edge pair connecting the roots
of $\tau_1$ and $\tau_2$), noting that $T(z)$ is the exponential generating function for rooted trees (since the number of rooted trees on $m$ labelled vertices is $m^{m-1}$),
and applying the multiplication formula for exponential generating functions (see for example~\cite[Prop.~5.1.1]{Sta99}).  Using~$|\Getree(m)|=(m-1)m^{m-2}=m^{m-1}-m^{m-2}=|\Gltree(m)|-|\Gtree(m)|$,
together with~\eqref{Fltree} and~\eqref{Fetree}, then gives~\eqref{Ftree}.

It is also a known result in graph theory (see for example~\cite[Eq.~(3.5)]{JanKnuLucPit93}), that
\begin{equation}\label{Qtree}\sum_{m=1}^\infty|\Gqtree(m)|\,\frac{z^m}{m!}=-\frac{T(z)}{2}-\frac{T(z)^2}{4}-\log\sqrt{1-T(z)},\end{equation}
where this can be obtained by first observing (as done similarly in Subsection~\ref{graphs})
that every quasitree corresponds uniquely to an undirected cycle $\tau_1,\ldots,\tau_\ell$ of disjoint rooted trees, for some $\ell\ge3$
(since a quasitree can be viewed as~$\tau_1,\ldots,\tau_\ell$, together with $\ell$ edges connecting the roots along the cycle).  Noting that the exponential
generating function for undirected cycles is
\begin{equation*}\sum_{\ell=3}^\infty\frac{(\ell-1)!}{2}\,\frac{z^\ell}{\ell!}=\frac{1}{2}\sum_{\ell=3}^\infty\frac{z^\ell}{\ell}=-\frac{1}{2}\Bigl(z+\frac{z^2}{2}+\log(1-z)\Bigr),\end{equation*}
and that the exponential generating function for rooted trees is $T(z)$, and applying the compositional formula for exponential generating functions (see for example~\cite[Prop.~5.1.4]{Sta99})
then gives~\eqref{Qtree}.

Substituting~\eqref{Fltree}--\eqref{Qtree} into~\eqref{fptz2} gives
\begin{align*}f_p(t,z)&=\frac{T(tz)}{t}-\frac{T(tz)^2}{2t}+p\,T(tz)+\frac{(t+1)\,T(tz)^2}{4t}-\frac{T(tz)}{2}-\frac{T(tz)^2}{4}-\log\sqrt{1-T(tz)}\\
&=\biggl(p-\frac{1}{2}+\frac{1}{t}\biggr)\,T(tz)-\frac{T(tz)^2}{4t}-\log\sqrt{1-T(tz)}.\end{align*}
Applying~\eqref{Fptz} then gives
\begin{equation*}F_p(t,z)=\frac{e^{(p-1/2+1/t)\,T(tz)-T(tz)^2/(4t)}}{\sqrt{1-T(tz)}},\end{equation*}
and using~\eqref{EhrFDef} with $p=n-m+1$ then gives
\begin{equation}\label{EhrT}
\ehr_{\P(m,n)}(t)\big|_{n\ge m-1}
=m!\,t^m\,[z^m]\,\frac{e^{(n-m+1/2+1/t)\,T(z)-T(z)^2/(4t)}}{\sqrt{1-T(z)}},\end{equation}
where, for a nonnegative integer~$k$ and power series $f(z)$, $[z^k]f(z)$ denotes the coefficient of~$z^k$ in the expansion of $f(z)$.

As noted in~\cite[Eq.~(2.38)]{CorGonHarJefKnu96}), for any nonnegative integer $k$ and power series $f(z)$, the function $T(z)$ satisfies
\begin{equation}\label{Tprop2}
[z^k]\,f\bigl(T(z)\bigr)=[z^k]\,f(z)\,(1-z)\,e^{kz}.\end{equation}
Applying~\eqref{Tprop2} to \eqref{EhrT}, with $f(z)=e^{(n-m+1/2+1/t)\,z-z^2/(4t)}\big/\sqrt{1-z}$, gives
\begin{equation}\label{EhrC}
\ehr_{\P(m,n)}(t)\big|_{n\ge m-1}=m!\,t^m\,[z^m]\,\sqrt{1-z}\,e^{(n+1/2+1/t)\,z-z^2/(4t)}.\end{equation}
Finally,~\eqref{EhrC} gives
\begin{equation*}
\ehr_{\P(m,n)}(t)\big|_{n\ge m-1}=m!\,t^m\sum_{0\le i\le j\le m}\bigl([z^i]\,\sqrt{1-z}\,\bigr)\,\bigl([z^{j-i}]\,e^{(n+1/2+1/t)\,z}\bigr)\,\bigl([z^{m-j}]\,e^{-z^2/(4t)}\bigr),
\end{equation*}
and explicitly evaluating the power series coefficients gives~\eqref{EhrEq}, after some straightforward simplification.

\section{Discussion}\label{disc}
Some comments on Theorem~\ref{EhrTh} are as follows.
\begin{list}{(\roman{c-disc})}{\usecounter{c-disc}\setlength{\labelwidth}{9.1mm}\setlength{\leftmargin}{9.1mm}\setlength{\labelsep}{1.5mm}\setlength{\topsep}{0.9mm}\setlength{\itemsep}{1mm}}
\item As shown in~\cite[Exer.~4.64 \& Sol.~4.64]{Sta12} and references within, the
Ehrhart polynomial of the permutohedron $\Pi(m,\ldots,1)$ can be expressed as an edge-weighted sum over all forests with vertices labelled $1,\ldots,m$.
Hence,~\eqref{EhrGra} provides an analogous expression for the partial permutohedron $\P(m,n)|_{n\ge m-1}$.  Also, an expression
for the Ehrhart polynomial of the permutohedron $\Pi(m,\ldots,1)$ in terms of the function $T(z)$ in~\eqref{Tdef} is
obtained in~\cite[1st Eq., Thm.~5.2]{ArdBecMcw20}, and~\eqref{EhrT} can be regarded as an analogous expression for the partial permutohedron $\P(m,n)|_{n\ge m-1}$.
\item As shown in~\cite[Rem.~6.6]{BehCasChaDiaEscHarInk22}, it follows from~\eqref{EhrC} that, for $n\ge m-1$,
\begin{multline}\label{recur}\qquad\ehr_{\P(m,n)}t)
=(mt+nt-t+1)\ehr_{\P(m-1,n)}(t)\\-(m-1)(nt+t/2+3/2)t\ehr_{\P(m-2,n)}(t)+(m-1)(m-2)t^2\ehr_{\P(m-3,n)}(t)/2,\end{multline}
where this recurrence determines $\ehr_{\P(m,n)}|_{n\ge m-1}$ with suitable initial conditions.
It would be interesting to find an alternative proof of Theorem~\ref{EhrTh} which obtains~\eqref{recur} directly.
\item As shown in~\cite[p.~36]{BehCasChaDiaEscHarInk22}, the explicit expression~\eqref{volP} for the volume of~$\P(m,n)_{n\ge m-1}$ follows
from~\eqref{EhrEq}.  Hence, Theorem~\ref{EhrTh} provides a new proof of~\eqref{volP}.
\item As discussed in~\cite[Rems.~6.4 \&~6.5]{BehCasChaDiaEscHarInk22}, Theorem~\ref{EhrTh} provides the answer
a question in~\cite{Sta20} for the number of integer points in the parking function polytope~$P_m$ (see~(iii) of Section~\ref{conn}), as well as to certain
other related questions, including one in~\cite[Section~6]{Sel22}.
\item It is expected that the $(a,b)$-parking function polytopes in~(iii) of Section~\ref{conn} can be expressed in terms of type-$\mathcal{Y}$ generalized permutohedra for arbitrary~$a$ and~$b$,
and that the techniques used to obtain Theorem~\ref{EhrTh} could also be used to obtain expressions for Ehrhart polynomials of such polytopes.
\end{list}

\let\oldurl\url
\makeatletter
\renewcommand*\url{%
        \begingroup
        \let\do\@makeother
        \dospecials
        \catcode`{1
        \catcode`}2
        \catcode`\ 10
        \url@aux
}
\newcommand*\url@aux[1]{%
        \setbox0\hbox{\oldurl{#1}}%
        \ifdim\wd0>\linewidth
                \strut
                \\
                \vbox{%
                        \hsize=\linewidth
                        \kern-\lineskip
                        \raggedright
                        \strut\oldurl{#1}%
                }%
        \else
                \hskip0pt plus\linewidth
                \penalty0
                \box0
        \fi
        \endgroup
}
\makeatother
\gdef\MRshorten#1 #2MRend{#1}
\gdef\MRfirsttwo#1#2{\if#1M
MR\else MR#1#2\fi}
\def\MRfix#1{\MRshorten\MRfirsttwo#1 MRend}
\renewcommand\MR[1]{\relax\ifhmode\unskip\spacefactor3000 \space\fi
  \MRhref{\MRfix{#1}}{{\tiny \MRfix{#1}}}}
\renewcommand{\MRhref}[2]{
 \href{http://www.ams.org/mathscinet-getitem?mr=#1}{#2}}

\bibliography{Bibliography}
\bibliographystyle{amsplainhyper}
\end{document}